\documentclass[journal]{IEEEtran}

\ifCLASSINFOpdf
\else
   \usepackage[dvips]{graphicx}
\fi

\usepackage{balance}

\usepackage[usenames]{color}
\usepackage[colorlinks = true]{hyperref}
\usepackage[english]{babel}
\usepackage{mathtools,algorithm,enumitem}
\usepackage{bm,bbm,amsmath,amssymb,amsthm,graphicx}
\usepackage{url}
\usepackage[caption=false,font=normalsize,labelfont=sf,textfont=sf]{subfig}
\usepackage{algpseudocode}
\usepackage{pifont}
\usepackage{todonotes}
\usepackage{tikz}
\usetikzlibrary{calc,patterns,decorations.pathmorphing,decorations.markings}

\usepackage{cleveref,cite}

\newcommand{\SLb}[1]{\textcolor{black}{#1}}      
\newcommand{\MBW}[1]{\textcolor{black}{#1}}      

\numberwithin{equation}{section}

\newcommand{\change}[1]{\ensuremath{\operatorname{#1}}}
\newcommand{\MAT}{\left[ \begin{array}}  
\newcommand{\mat}{\end{array} \right]}
\newtheorem{Definition}{Definition}[section]
\newtheorem{Corollary}{Corollary}[section]
\newtheorem{Lemma}{Lemma}[section]
\newtheorem{Theorem}{Theorem}[section]
\newtheorem{Example}{Example}[section]
\newtheorem{Q}{Question}[section]

\def \st {\operatorname*{s.t. }}

\def \b {\bm{b}}

\def \D {\mathbf{D}}

\def \Dt {\widetilde{\mathbf{D}}}
\def \e {\bm{e}}

\def \G {\mathbf{G}}

\def \I {\mathbf{I}}

\def \NN {\mathcal{N}}

\def \Pbf {\mathbf{P}}

\def \Q {\mathbf{Q}}

\def \rs{r^\star}
\def \R {\mathbf{R}}

\def \RRR {\mathbb{R}}

\def \U {\mathbf{U}}

\def \Ut {\widetilde{\mathbf{U}}}
\def \Us {\mathbf{U}^\star}

\def \V {\mathbf{V}}

\def \Vt {\widetilde{\mathbf{V}}}

\def \X {\mathbf{X}}

\def \Xs {\mathbf{X}^\star}

\def \sumj{\sum_{j=1}^J}

\def \bSigma {\boldsymbol{\Sigma}}

\def \zero {\mathbf{0}}

\newcommand{\tmop}[1]{\ensuremath{\operatorname{#1}}}

\newcommand{\vct}[1]{\mathbf{#1}}
\newcommand{\mtx}[1]{\mathbf{#1}}

\newcommand{\set}[1]{\mathcal{#1}}

\def \mPi     {{\boldsymbol{\Pi}}}

\def \b {\vct{b}}

\def \e {\vct{e}}

\def \mD {\mtx{D}}

\def \mG {\mtx{G}}

\def \mU {\mtx{U}}
\def \mV {\mtx{V}}

\def \mX {\mtx{X}}

\def \calA {\set{A}}

\def \calG {\set{G}}

\DeclareMathOperator*{\minimize}{\operatorname{minimize}}

\def \ST { \;\; \mbox{subject to} \;\; }
\def \nn       {\nonumber}

\def \zero     {\mathbf{0}}

\def \lg       {\left\langle}
\def \rg       {\right\rangle}

\def \rank     {\tmop{rank}}

\def \Null     {\tmop{null}}

\def \det      {\tmop{det}}

\def \R     {\mathbb{R}}

\begin{document}

\title{The Global Geometry of Centralized and Distributed Low-rank Matrix Recovery without Regularization}

\author{Shuang Li, Qiuwei Li, Zhihui Zhu, Gongguo Tang, and Michael B. Wakin
\thanks{S. Li, Q. Li, G. Tang and M. B. Wakin are with the Department of Electrical Engineering, Colorado School of Mines, Golden, CO 80401, USA (e-mail: \{shuangli, qiuli, gtang, mwakin\}@mines.edu). Z. Zhu is with the Ritchie School of Engineering and Computer Science, University of Denver, Denver, CO 80208, USA (e-mail: zhihui.zhu@du.edu). This work was supported by NSF grant CCF-1704204, and the DARPA Lagrange Program under ONR/SPAWAR contract N660011824020.}
}



\maketitle

\begin{abstract}
Low-rank matrix recovery is a fundamental problem in signal processing and machine learning. A recent very popular approach to recovering a low-rank matrix $\X$ is to factorize it as a product of two smaller matrices, i.e.,  $\X = \mU\mV^\top$, and then optimize over $\U, \V$ instead of $\mX$. Despite the resulting non-convexity, recent results have shown that \MBW{many factorized objective functions actually} have benign global geometry---with no spurious local minima and satisfying the so-called strict saddle property---ensuring convergence to a global minimum for many local-search algorithms. \MBW{Such results hold whenever the original objective function is restricted strongly convex and smooth.} However, most of these results actually consider a modified cost function that includes a balancing regularizer. While useful for deriving theory, this balancing regularizer does not appear to be necessary in practice. In this work, we close this theory-practice gap by proving that the \MBW{unaltered} factorized non-convex problem, without the balancing regularizer, also has similar benign global geometry. Moreover, we also extend our theoretical results to the field of distributed optimization.


\end{abstract}


\IEEEpeerreviewmaketitle

\section{Introduction}

\IEEEPARstart{I}{n} the problem of low-rank matrix recovery, a great number of efforts have been made to minimize a 
loss function $f(\mX)$ over the non-convex rank constraint $\rank(\mX)\le r$, where $\X\in\RRR^{n\times m}$ and $r \ll \min\{n,m\}$. Among which, a popular way is to replace the rank constraint with the Burer-Monteiro factorization, i.e., $\X = \U\V^\top$ with $\U\in\RRR^{n\times r}$ and $\V\in\RRR^{m\times r}$~\cite{burer2003nonlinear,burer2005local}, changing the objective function from $f(\mX)$ to $g(\mU,\mV) = f(\mU\mV^\top)$. This factorization approach can often lead to lower computational and storage complexity, while raising new questions about whether an algorithm can converge to favorable solutions since the bilinear form $\U\V^\top$ naturally introduces non-convexity.
Fortunately, it is observed that simple iterative algorithms find global optimal solutions in many low-rank matrix recovery problems \cite{li2019bregman,chi2019nonconvex,li2018non,li2019landscape, LiEtAl2017Geometry,zhu2017global,ge2017no,zhu2018global,zhu2018global2,park2017non}.

Recent years have seen a surge of interest in understanding these surprising phenomena by analyzing the landscape of the factorized cost function $g(\mU,\mV)$. To accomplish this, many existing works \cite{ge2017no,zhu2017global,zhu2018global,zhu2018global2,park2017non,wang2017unified,tu2016low,zheng2016convergence,park2016finding} actually add a balancing regularizer
\begin{align}
\tmop{R}(\mU,\mV)\doteq\|\mU^\top\mU-\mV^\top\mV\|_F^2,
\label{eq:balance_reg}
\end{align}
which implicitly forces $\U$ and $\V$ to have equal energy,
to the objective function $g(\U,\V)$. These works then show \MBW{that, for broad classes of problems,} the regularized cost functions have a benign geometry,  where every local minimum is a global minimum and every first-order critical point is either a local minimum or a strict saddle~\cite{sun2015complete,ge2015escaping}. This favorable property ensures a convergence to a global minimum for many local search methods \cite{ge2015escaping,lee2016gradient,jin2017escape,li2019alternating,nesterov2006cubic,more1983computing,lu2019pa}.

\subsection{What Is The Role of The Balancing Regularizer?}

If $(\U,\V)$ is a critical point of $g(\U,\V)$, then $({\mU}\mG,\mV\mG^{-\top})$ is also a critical point for any invertible $\G\in\R^{r\times r}$. This scaling ambiguity in the critical points can result in an infinite number of connected critical points including those ill-conditioned points when $\|\mG\|_F$ goes to 0 or $\infty$, which could bring new challenges in analyzing the geometric landscape as one must analyze the optimality of any critical point. In order to remove this ambiguity, many researchers \cite{ge2017no,zhu2017global,zhu2018global,zhu2018global2,park2017non,wang2017unified,tu2016low,zheng2016convergence,park2016finding} utilize the balancing regularizer~\eqref{eq:balance_reg}. In particular, it has been shown that adding the regularizer~\eqref{eq:balance_reg} forces all critical points $(\mU,\V)$ to be balanced, i.e., $\mU^\top\mU = \mV^\top\mV$.


\subsection{Is The Balancing Regularizer Really Necessary?}
Most previous works add the balancing regularizer~\eqref{eq:balance_reg} to the cost function in order to simplify the landscape analysis. However, we have observed that one can achieve almost the same performance without adding the balancing regularizer in \cite{zhu2018global}. Also, in practice this additional regularizer is rarely utilized~\cite{mabeyond2019}, which implies a gap between theory and practice.
This naturally raises the main question that will be addressed in this work:  Is the balancing regularizer~\eqref{eq:balance_reg} truly necessary? In other words, can we characterize the global geometry of the factorization approach without the balancing regularizer?

Several works \cite{du2018algorithmic,NIPS2019_8960,mabeyond2019,Ma2019} answer this question by analyzing the behavior of gradient descent on some particular optimization problems, and show that the iterates of gradient descent stay in the (approximately) balanced path from some specific initialization and finally converge to a global optimal solution. However, these results are restricted to gradient descent with a specific initialization.
There are also some works that analyze the geometric landscape of some specific optimization problems, such as matrix factorization~\cite{zhu2019distributed}, or linear neural network optimization~\cite{nouiehed2018learning,zhu2018global2,kawaguchi2016deep}.

In this work, we answer this question by directly analyzing the landscape 
\MBW{of the unaltered factorized non-convex problem, {\em without} the balancing regularizer~\eqref{eq:balance_reg}. In particular, over the general class of problems where the cost function $f$ is restricted strongly convex and smooth (see \Cref{def:rscs}), we show under mild conditions that any critical point of the factorized cost function $g$ (including any unbalanced critical point) is either a global optimum or a strict saddle. This helps close the theory-practice gap and resolves the open problem in \cite{zhu2018global}. Moreover, we extend our results to the corresponding distributed setting and show that many global consensus problems inherit the benign geometry of their original centralized counterpart.}

\MBW{Before proceeding,} we present a toy example to illustrate our main observation.

\begin{figure}[t]
\begin{center}
\hspace*{-.4cm}
\makebox{
\begin{minipage}{0.13\textwidth}
\includegraphics[width=.9\textwidth]{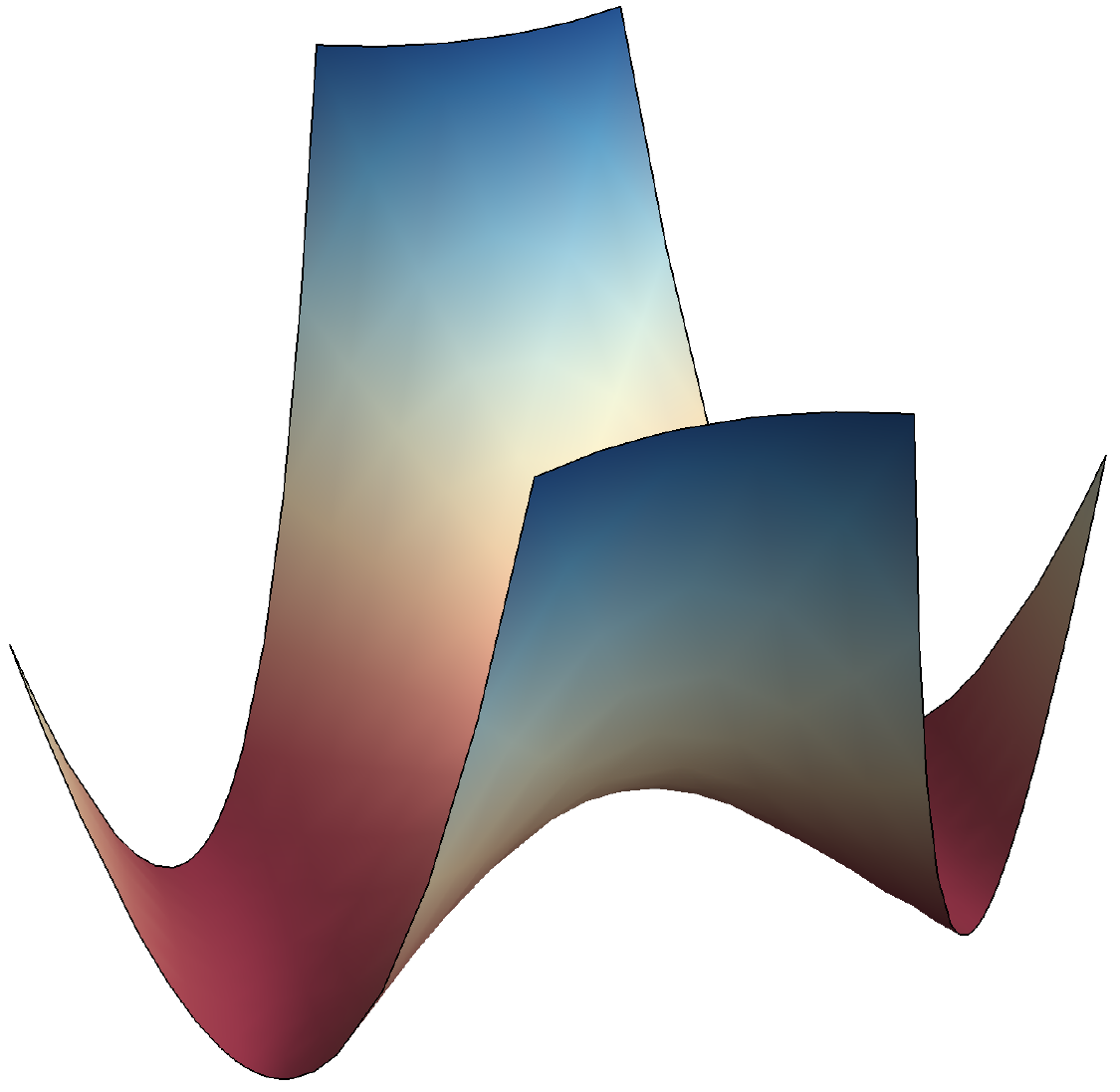}
\end{minipage}
\!\!\!\!\!\!\!\!
\begin{minipage}{0.13\textwidth}
\includegraphics[width=\textwidth]{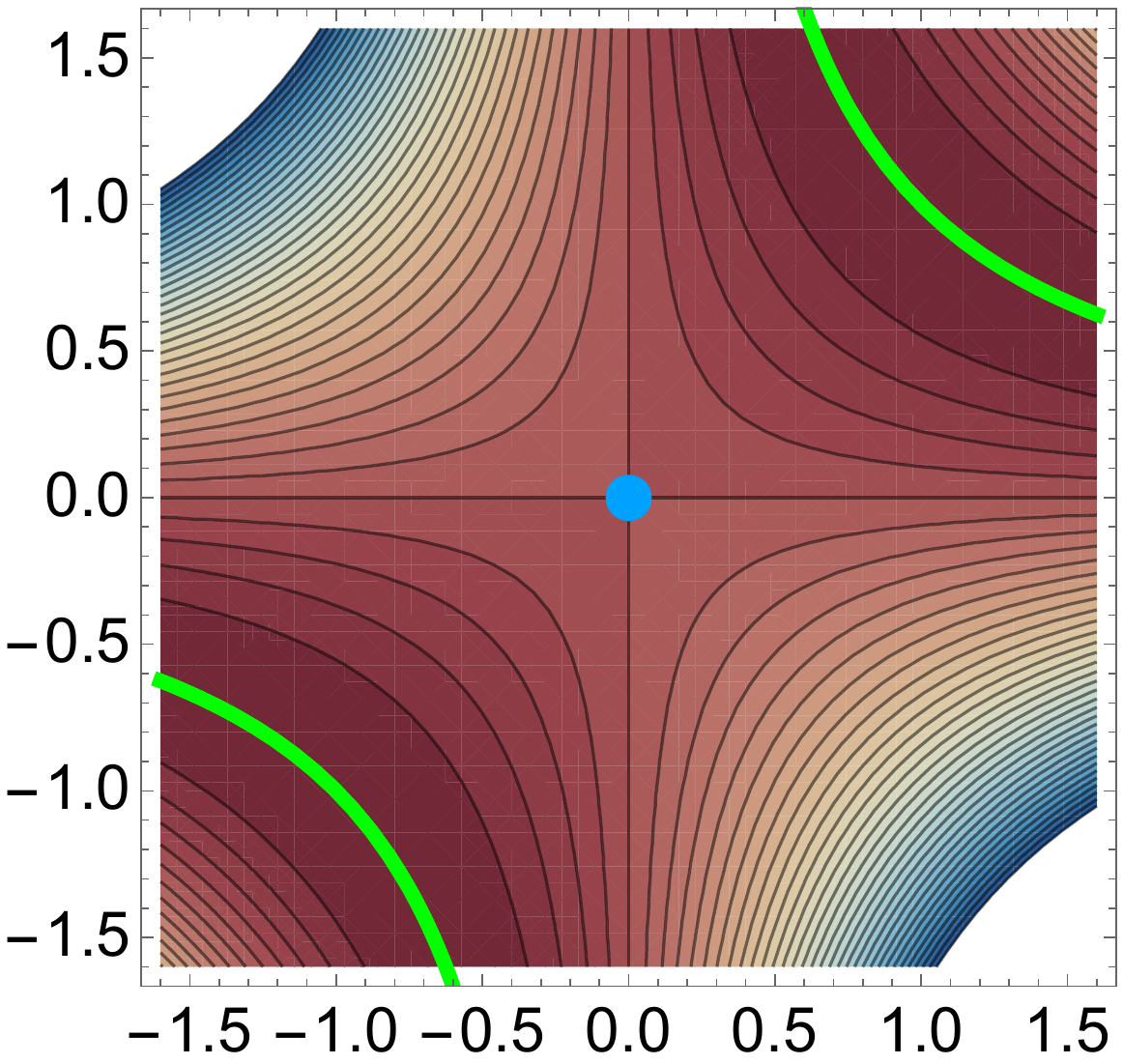}
\end{minipage}
\!\!\!\!\!
\begin{minipage}{0.13\textwidth}
\centering
\includegraphics[width=.9\textwidth]{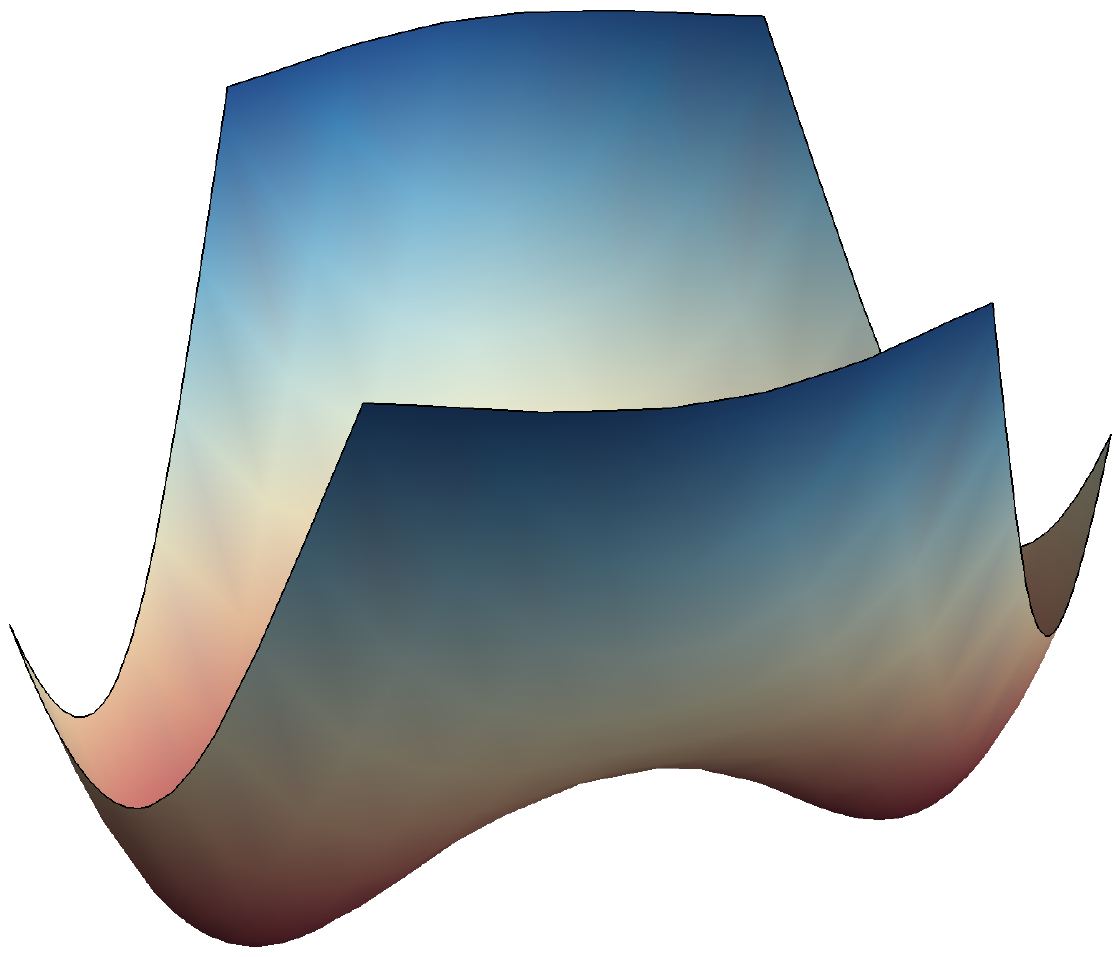}
\end{minipage}
\!\!\!\!\!
\begin{minipage}{0.13\textwidth}
\centering
\includegraphics[width=\textwidth]{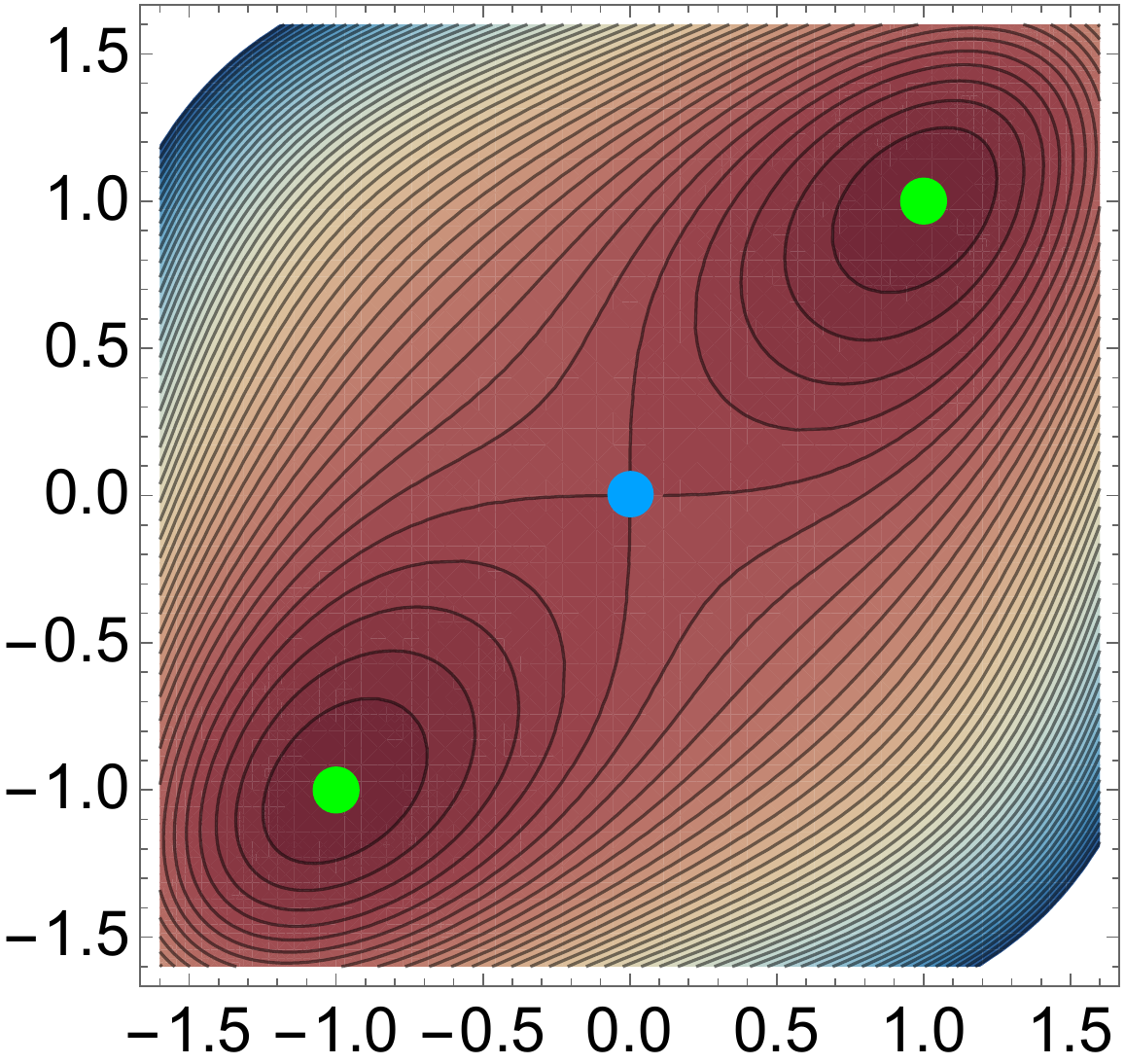}
\end{minipage}}
\\[.1cm]
\hspace*{-.7cm}
\makebox{
\begin{minipage}{0.13\textwidth}
\centering{\footnotesize{(a)}}
\end{minipage}
\begin{minipage}{0.13\textwidth}
\centering{\footnotesize{(b)}}
\end{minipage}
\begin{minipage}{0.13\textwidth}
\centering
\centering{\footnotesize{(c)}}
\end{minipage}
\begin{minipage}{0.13\textwidth}
\centering
\centering{\footnotesize{(d)}}
\end{minipage}}
\end{center}
\vspace{-.4cm}
\caption{(a, b) and (c, d) are the landscapes of the non-regularized function $g(u,v)$ and the regularized cost function $\widetilde{g}(u,v)$ with $\mu=\frac{1}{16}$, respectively.
One can observe that although $g$ has an infinite number of (connected) critical points while $\widetilde{g}$ has just three critical points at $\pm(1,1)$ and $(0,0)$, both cost functions have benign landscapes since any critical point is either a global minimizer or a strict saddle. The points marked with green and blue in (b, d) denote the global minimizers and saddle points, respectively.}
\label{fig:toyexample}
\end{figure}

\begin{Example}[Matrix factorization -- the scalar case]
Consider
an asymmetric matrix factorization cost function
$g(u,v)\doteq \frac{1}{2}(1- u v)^2$,
whose critical points $(u,v)$ satisfy $uv=1$ or $(u,v) = (0,0)$. The critical points of the corresponding regularized function
$\widetilde{g}(u,v)\doteq \frac{1}{2}(1- u v)^2+\frac{\mu}{4}(u^2-v^2)^2$
with some $\mu > 0$ satisfy $(uv-1)v+\mu(u^2-v^2)u=0$ and $(uv-1)u-\mu(u^2-v^2)v=0$,
which gives only three critical points $(1,1),(-1,-1)$ and $(0,0)$. Therefore, for $\widetilde{g}(u,v)$, one only needs to check the Hessian evaluated at these three critical points:
$
\nabla^2\widetilde{g}(1,1)=\nabla^2\widetilde{g}(-1,-1)=\left[\begin{smallmatrix}
1+2\mu& 1-2\mu \\ 1-2\mu & 1+2\mu
\end{smallmatrix}\right]
\succ 0,
$
and
$
\nabla^2\widetilde{g}(0,0)=\left[\begin{smallmatrix*}[r]
0&-1\\-1 &0
\end{smallmatrix*}\right]	
$
which has a strictly negative eigenvalue $-1$. Thus any critical point of $\widetilde{g}(u,v)$ is either a global minimum or a strict saddle, which implies a favorable landscape of the regularized cost function $\widetilde{g}(u,v)$. As can be seen, adding the balancing regularizer can largely simplify the landscape analysis. However, this does not imply that the original function $g(u,v)$ does not have a benign geometry. Indeed, one can observe that any critical point of $g(u,v)$ either satisfies $uv=1$ (globally optimal) or $(u,v)=(0,0)$ (strict saddle since $\nabla^2 g(0,0)=\left[\begin{smallmatrix*}[r]
0&-1\\-1 &0
\end{smallmatrix*}\right]$ has a negative eigenvalue $-1$). The landscapes of $g$ and $\widetilde{g}$ are shown in  \Cref{fig:toyexample}.
\end{Example}




The remainder of this paper is organized as follows. In Section~\ref{sec:prob}, we formulate the problems in both centralized and distributed settings. We present our main theorem and its proof in Section~\ref{sec:main}. In Section~\ref{sec:simu}, we conduct a series of experiments to further support our theory. Finally, we conclude our work in Section~\ref{sec:conc}.

\section{Problem Formulation}
\label{sec:prob}

We first consider the following problem of minimizing a general objective function over the set of low-rank matrices:
\begin{align}
\minimize_{\mX\in\R^{n\times m}}~f(\mX) \quad \ST \rank(\mX)\le r,
\label{problem:with:rank}
\end{align}
which is a fundamental problem that often appears in the fields of signal processing and machine learning. Plugging the Burer-Monteiro type decomposition~\cite{burer2003nonlinear,burer2005local}, i.e., $\X = \U\V^\top$ with $\U\in\RRR^{n\times r}$ and $\V\in\RRR^{m\times r}$, into the above cost function, one can remove the low-rank constraint and get the following unconstrained optimization
\begin{align}
\minimize_{\U\in\R^{n\times r},\V\in\R^{m\times r}}~g(\U,\V)  \doteq  f(\U\V^\top),
\label{eq:cost_nonconv}	
\end{align}
which is a non-convex optimization problem \MBW{we refer to} as {\em centralized} low-rank matrix recovery.  The above optimization appears in many applications including low-rank matrix approximation~\cite{zhu2017global}, matrix sensing \cite{ge2017no}, matrix completion \cite{sun2016guaranteed}, and linear neural network optimization~\cite{nouiehed2018learning,zhu2018global2,kawaguchi2016deep}.
Note that in centralized low-rank matrix recovery, all the computations happen at one ``central'' node that has full access, for example, to the data matrix or the measurements.

In the second part of this work, we study the impact of {\em distributing} the centralized low-rank matrix recovery problem for general cost functions.  Consider a separable cost function
$f(\U\V^\top) = \sumj   f_j(\mU\mV_j^\top)$,
 where $\mU\in\R^{n\times r}$ is the common variable in all of the objective functions $\{f_j\}_{j\in[J]}$ and $\mV_j\in\RRR^{m_j \times r }$ as a submatrix of $\mV=[\mV_1^\top~\cdots~\mV_J^\top]^\top\in\R^{m\times r}$ is the local variable only corresponding to objective function $f_j$.
Then, the centralized optimization~\eqref{eq:cost_nonconv}  becomes
\begin{align}
\minimize_{\U,\{\mV_j\}}~
\sumj   f_j(\mU\mV_j^\top).
\label{eq:cost_nonconv_dis00}	
\end{align}
In the {\em distributed} setting, one distributes \eqref{eq:cost_nonconv_dis00} across a network of $J$ agents and considers the following optimization
\begin{align}
\minimize_{\{\mU^j\},\{\mV_j\}}  \sumj f_j(\mU^j\mV_j^\top)~\st~\mU^1=\cdots=\mU^J.
\label{ECD}
\end{align}
Here, $\mU^j$ and $\mV_j$ are the so-called consensus and local variables at node $j$.
In this work, we consider the above equality-constrained distributed problem \eqref{ECD} by reformulating it as the following unconstrained optimization problem
\begin{align}
\minimize_{\{\mU^j\},\{\mV_j\}}  \sumj f_j(\mU^j\mV_j^\top)\!+\!\!\!\!\!\sum_{(i,j)\in\calG}\!\!\!\!\! w_{i,j} \|\mU^j\!-\!\mU^i\|_F^2.
\label{eq:cost_nonconv_dis0:ECD:Penalty}	
\end{align}
 Here, $\calG$ denotes any connected network over $[J]^2$ with \SLb{$[J] \doteq \{1,\cdots,J\}$ and $[J]^2 = [J]\times [J]$}~\cite{zhu2019distributed}, and $\{w_{i,j}\}_{(i,j)\in\calG}$ are symmetric positive weights, i.e., $w_{j,i}=w_{i,j}> 0$.   The second term is added to the objective function for the purpose of promoting equality among the consensus variables $\mU^j$.

In this work, our main goal is to characterize the global geometry of the non-convex {\em centralized} cost function~\eqref{eq:cost_nonconv} and non-convex {\em distributed} cost function~\eqref{eq:cost_nonconv_dis0:ECD:Penalty}. In particular, we show that under the same assumptions as  required in the previous works, any critical point is either a global minimum or a strict saddle, where the Hessian has a strictly negative eigenvalue, without adding the balancing regularizer.


\section{Main Results}
\label{sec:main}

\subsection{Landscape of Centralized Low-rank Matrix Recovery}
\label{sec:land_cen}

In this subsection, we present the geometric landscape of the {\em centralized} optimization~\eqref{eq:cost_nonconv}. We start by introducing the restricted strongly convex and smooth property.


\begin{Definition} \label{def:rscs}
(\cite{zhu2018global,li2018non})
A function $f(\X)$ is said to be $(2r,4r)$-restricted strongly convex and smooth if
\begin{align}
\alpha \|\D\|_F^2 \leq \nabla^2 f(\X)[\D,\D] \leq \beta \|\D\|_F^2	
\label{eq:RSC}\end{align}
holds for any matrix $\X\in\RRR^{n\times m}$ with rank at most $2r$ and $\D\in\RRR^{n\times m}$ with rank at most $4r$. Here, $\alpha$ and $\beta$	are some positive constants, \SLb{ and $ \nabla^2 f(\X)[\D,\D] = \sum_{i,j,k,l} \frac{\partial^2 f(\X)}{\partial \X_{ij} \partial \X_{kl}} \D_{ij} \D_{kl}$ denotes a bilinear form of the Hessian of $f(\X)$.}
\end{Definition}

\MBW{Unlike the standard strongly convex and smooth condition which requires \eqref{eq:RSC} to hold for any $\mX$ and $\mD$, the above restricted version only requires \eqref{eq:RSC} to hold for low-rank matrices, making it amenable for low-rank matrix recovery problems. For example, in matrix sensing the goal is to recover a low-rank matrix $\mX$ from linear measurements $\calA(\mX)$. The linear operator $\calA$ often satisfies the restricted isometry property (RIP), which can be interpreted as satisfying~\eqref{eq:RSC} for all low-rank matrices $\D$ and all $\X$; see \cite{zhu2018global} for details.}
\begin{Theorem} \label{thm:local_mtxsen}

Assume that the cost function $f(\X)$ in~\eqref{problem:with:rank} satisfies the $(2r,4r)$-restricted strongly convex and smooth property with positive constants $\alpha$ and $\beta$ satisfying $\beta/\alpha\le 3/2$. Also assume that $f(\X)$ has a critical point $\mX^\star$ with $\rank(\mX^\star)\le r$. Then, any critical point $(\mU,\mV)$ of $g(\mU,\mV)$ in~\eqref{eq:cost_nonconv} is either a global minimum (i.e., $\U\V^\top = \Xs$) or a strict saddle (i.e., $\lambda_{\min}(\nabla^2 g(\U,\V))<0$).
\end{Theorem}

\vspace{-.05in}
\begin{proof}
It follows from~\cite[Proposition 1]{zhu2018global} that the critical point $\Xs$ of $f(\X)$ with $\change{rank}(\Xs) = \rs \leq r$ is its global minimum, namely, $f(\Xs) \leq f(\X)$ holds for any $\X\in\RRR^{n\times m}$ with $\change{rank}(\X)\leq r$. Moreover, the equality holds only at $\X = \Xs$. Then, for any critical point $(\U,\V)$ with $\U\V^\top = \Xs$, we have
$g(\U,\V) = f(\Xs),$ and hence $(\U,\V)$ is a global minimum.

For any critical point $(\U,\V)$ with $\U\V^\top \neq \Xs$, we next show that there exists a direction $\D\in\RRR^{(n+m)\times r}$ such that $\nabla^2 g(\U,\V) [\D,\D] < 0$, namely, $(\U,\V)$ is a strict saddle of $g(\U,\V)$. The remaining part of this proof is inspired by the proof of~\cite[Lemma 11.3]{zhu2019distributed} and~\cite[Theorem 8]{nouiehed2018learning} and is split into  two cases: 1)  $\rank(\U\V^\top) = r$, and 2)  $\rank(\U\V^\top) < r$.


\vspace{0.15cm}
\noindent{\bf Non-degenerate case: $\rank(\U\V^\top) = r$}

\vspace{0.15cm}

Let $\U\V^\top = \Pbf \bSigma \Q^\top$ be an SVD of $\U\V^\top$. It follows from $\change{rank}(\U\V^\top) = r$ that $\change{rank}(\U) = \change{rank}(\V) = r$, which further implies that $\U^\top\U$ and $\V^\top\V$ are invertible. Then, we define two matrices
$\G_1 \doteq (\U^\top \U)^{-1} \U^\top \Pbf \bSigma^{1/2}$, and $\G_2 \doteq (\V^\top \V)^{-1} \V^\top \Q \bSigma^{1/2}$.
It can be seen that $\G_1 \G_2^\top = \I_r$. 
\textcolor{black}{We also define balanced factors
$\Ut \doteq  \Pbf \bSigma^{1/2}$ and	$ \Vt \doteq  \Q\bSigma^{1/2}$.}

It can be seen that the new matrix pair $(\Ut,\Vt)$ satisfies
\begin{align}
\Ut \Vt^\top &= \U \V^\top, &
\Ut^\top \Ut &= \Vt^\top \Vt.	
\label{property:tilde}
\end{align}
Recall that for any critical point $(\U,\V)$ of $g = f(\U\V^\top)$, we have $\nabla g(\U,\V)=\zero$, i.e.,
$
\nabla f(\U\V^\top) \V = \zero$ and
$(\nabla f(\U\V^\top))^\top \U = \zero$.	
Together with the equalities in~\eqref{property:tilde}, we get
$\nabla \left(g+\frac{\mu}{4}\tmop{R}\right)\!(\widetilde{\mU},\widetilde{\mV})
\!=\!
\begin{bmatrix}
\begin{smallmatrix}
\nabla f(\widetilde{\mU}\widetilde{\mV}^\top)\widetilde{\mV}+\mu\widetilde{\mU}(\widetilde{\mU}^\top\widetilde{\mU}-\widetilde{\mV}^\top\Vt)
\\
(\nabla f(\widetilde{\mU}\widetilde{\mV}^\top))^\top\widetilde{\mU}-\mu \widetilde{\mV}(\widetilde{\mU}^\top\widetilde{\mU}-\widetilde{\mV}^\top\widetilde{\mV})	
\end{smallmatrix}
\end{bmatrix}
\!=\! \zero,$
where $\tmop{R}(\cdot)$ is the balancing regularizer introduced in~\eqref{eq:balance_reg} and $\mu>0$ is a regularizer parameter. This immediately implies that the new matrix pair $(\Ut,\Vt)$ is a critical point of the regularized cost function $g(\U,\V) + \frac{\mu}{4} \tmop{R}(\U,\V)$.

On the other hand, it follows from \cite{zhu2018global} that there exists a matrix $\Dt = [\Dt_{\Ut}^\top~\Dt_{\Vt}^\top]^\top \in\RRR^{(n+m)\times r}$ such that
\begin{align}
\nabla^2 \left(g(\widetilde{\mU},\widetilde{\mV})+\frac{\mu}{4}\tmop{R}(\widetilde{\mU},\widetilde{\mV})\right) \left[\widetilde{\mD},\widetilde{\mD}\right]<0
\label{negative:curvature}
\end{align}
holds for any $\Ut{\Vt}^\top \neq \Xs$.

Construct
$
\D = \left[\begin{smallmatrix} \D_{\U}\\ \D_{\V}  \end{smallmatrix} \right] = \left[\begin{smallmatrix} \Dt_{\Ut}\G_1^{-1}\\ \Dt_{\Vt}\G_2^{-1}  \end{smallmatrix} \right]$,
and denote $\mPi\doteq \mU\mD_\mV^\top+\mD_\mU\mV^\top$ and $\widetilde{\mPi}\doteq \widetilde{\mU}\widetilde{\mD}_{\Vt}^\top+\widetilde{\mD}_{\Ut}\widetilde{\mV}^\top$. Note that $\mPi = \widetilde{\mPi}$. It follows from~\eqref{negative:curvature} that
\begin{align*}
0>&\nabla^2 \left(g(\widetilde{\mU},\widetilde{\mV})+({\mu}/{4})\tmop{R}(\widetilde{\mU},\widetilde{\mV})\right)\left[\widetilde{\mD},\widetilde{\mD}\right]
\nn\\
=&
\lg 2\nabla f(\widetilde{\mU}\widetilde{\mV}^\top), \widetilde{\mD}_{\Ut}\widetilde{\mD}_{\Vt}^\top\rg +  \nabla^2f(\widetilde{\mU}\widetilde{\mV}^\top)\left[\widetilde{\mPi},\widetilde{\mPi}\right]	
\nn\\
&+ (\mu/2) \|\Dt_{\Ut}^\top\Ut+\Ut^\top\Dt_{\Ut}-\Dt_{\Vt}^\top\Vt-\Vt^\top\Dt_{\Vt}\|_F^2
\nn\\
\ge& \lg 2\nabla f(\widetilde{\mU}\widetilde{\mV}^\top), \widetilde{\mD}_{\Ut}\widetilde{\mD}_{\Vt}^\top\rg +  \nabla^2f(\widetilde{\mU}\widetilde{\mV}^\top)\left[\widetilde{\mPi},\widetilde{\mPi}\right]	
\nn\\
=& \lg 2\nabla f(\mU\mV^\top), \mD_\mU\mD_\mV^\top\rg+  \nabla^2f(\mU\mV^\top)\left[\mPi,\mPi\right]
\\
=& \nabla^2 g(\mU,\mV)[\mD,\mD],
\end{align*}
which further implies that any non-degenerate critical point $(\mU,\mV)$ with $\mU\mV^\top\neq\mX^\star$ is a strict saddle.

\begin{figure*}[t]
\begin{minipage}{0.19\linewidth}
\centering
\includegraphics[width=1.4in]{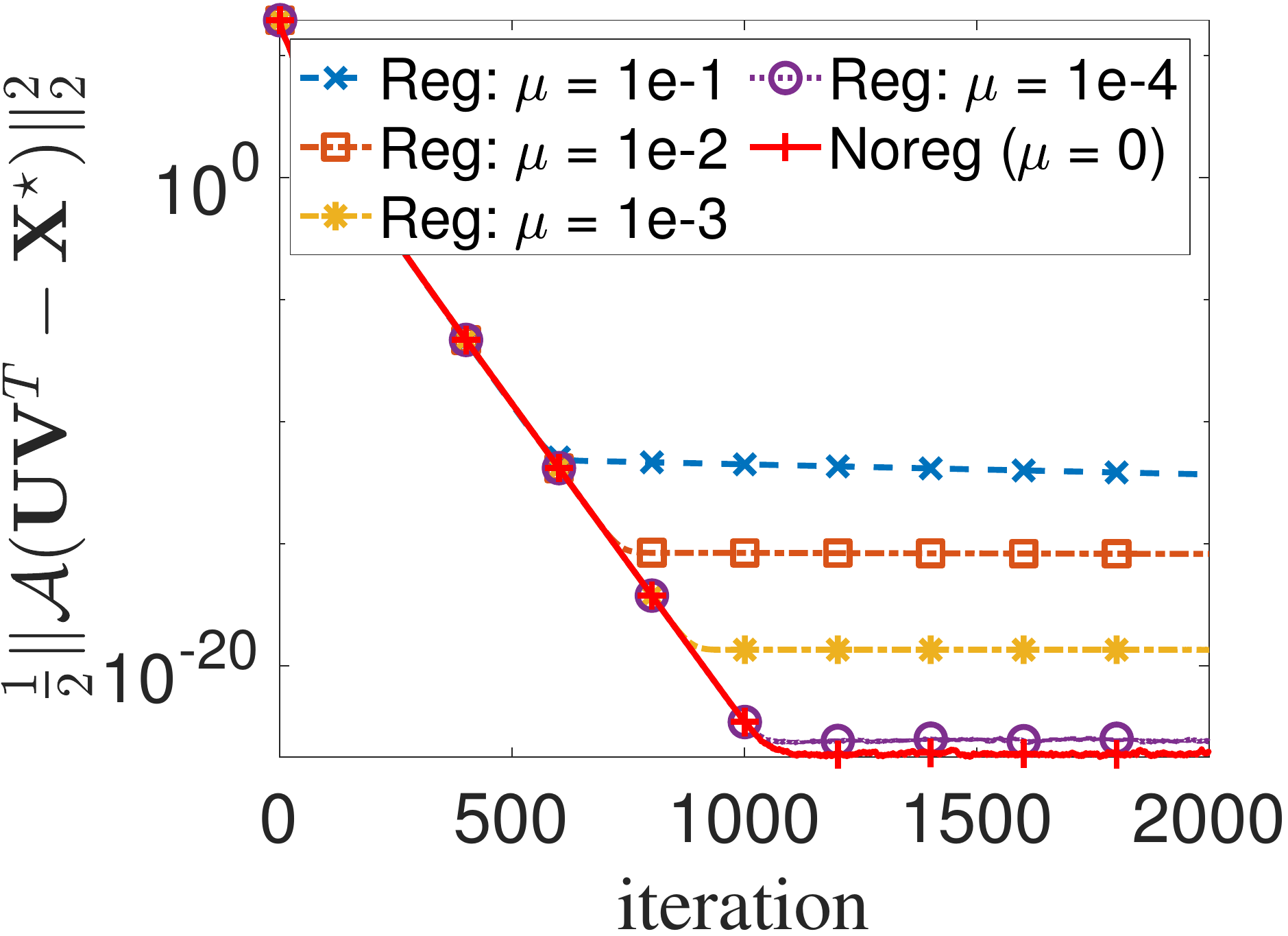}
\centerline{\footnotesize {(a) Centralized}}
\end{minipage}
\hfill
\begin{minipage}{0.19\linewidth}
\centering
\includegraphics[width=1.4in]{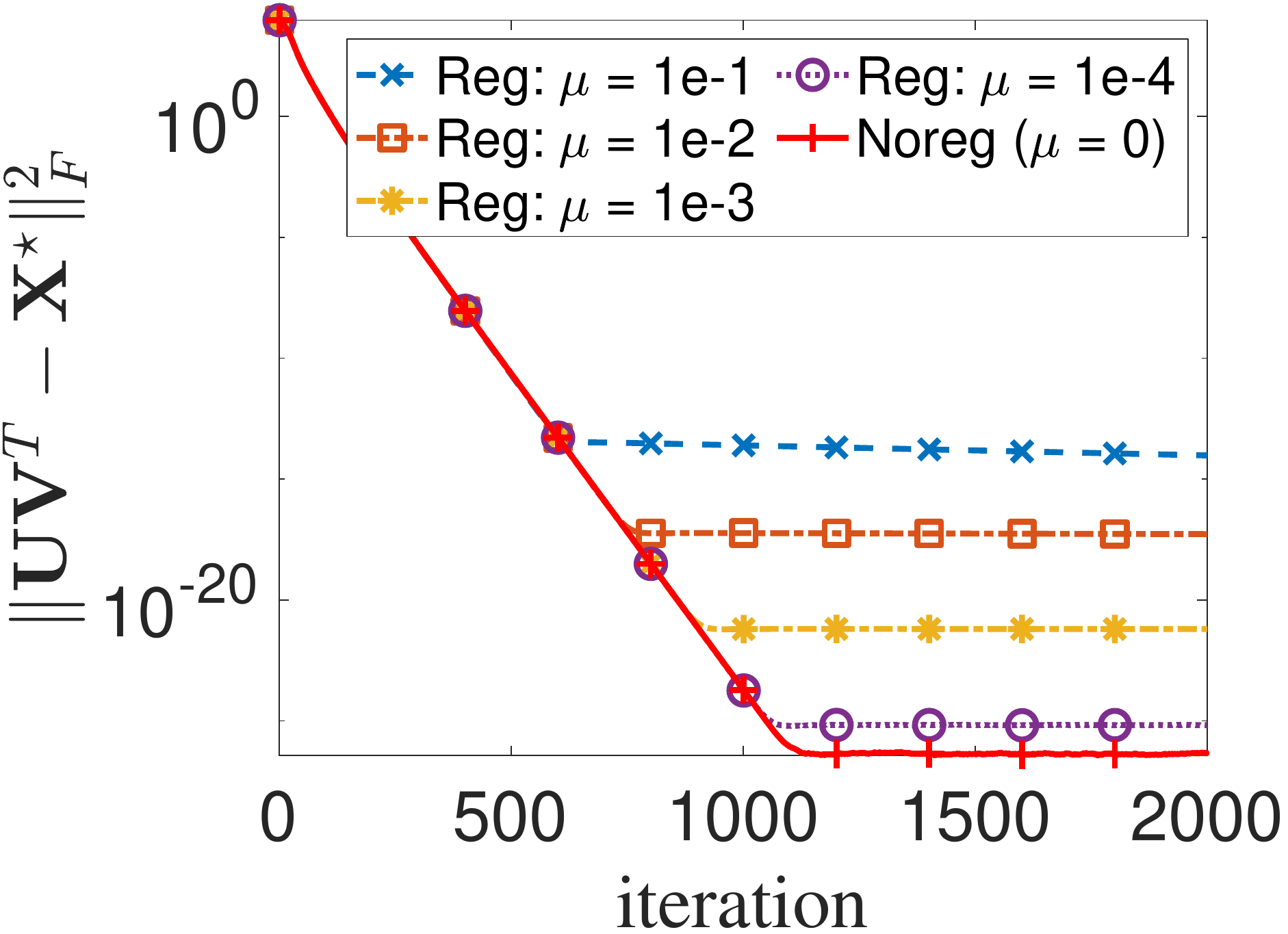}
\centerline{\footnotesize{(b) Centralized}}
\end{minipage} 
\hfill 
\begin{minipage}{0.19\linewidth}
\centering
\includegraphics[width=1.4in]{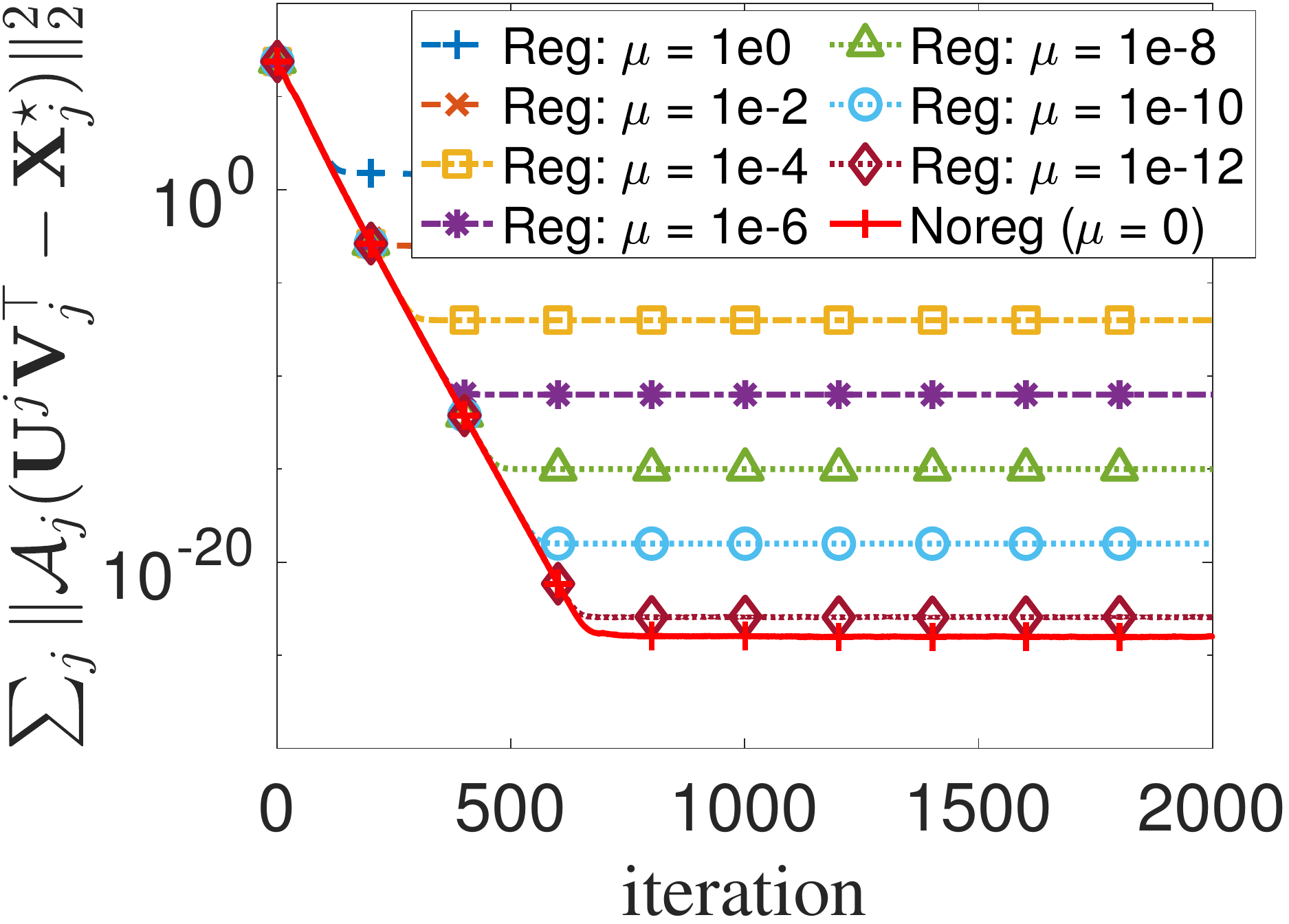}
\centerline{\footnotesize{(c) Distributed}}
\end{minipage}
\hfill
\begin{minipage}{0.19\linewidth}
\centering
\includegraphics[width=1.4in]{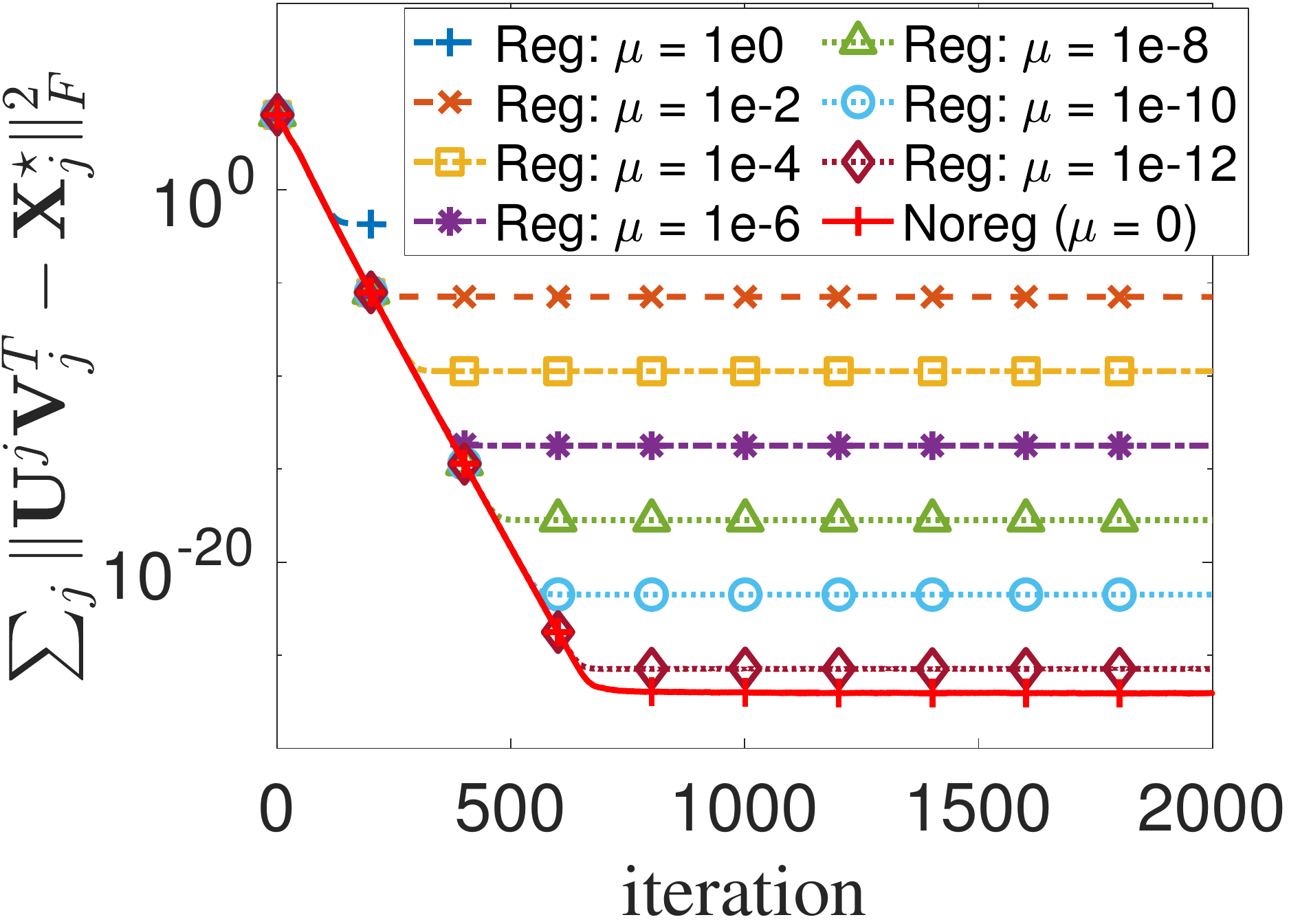}
\centerline{\footnotesize{(d) Distributed}}
\end{minipage}
\hfill
\begin{minipage}{0.19\linewidth}
\centering
\includegraphics[width=1.4in]{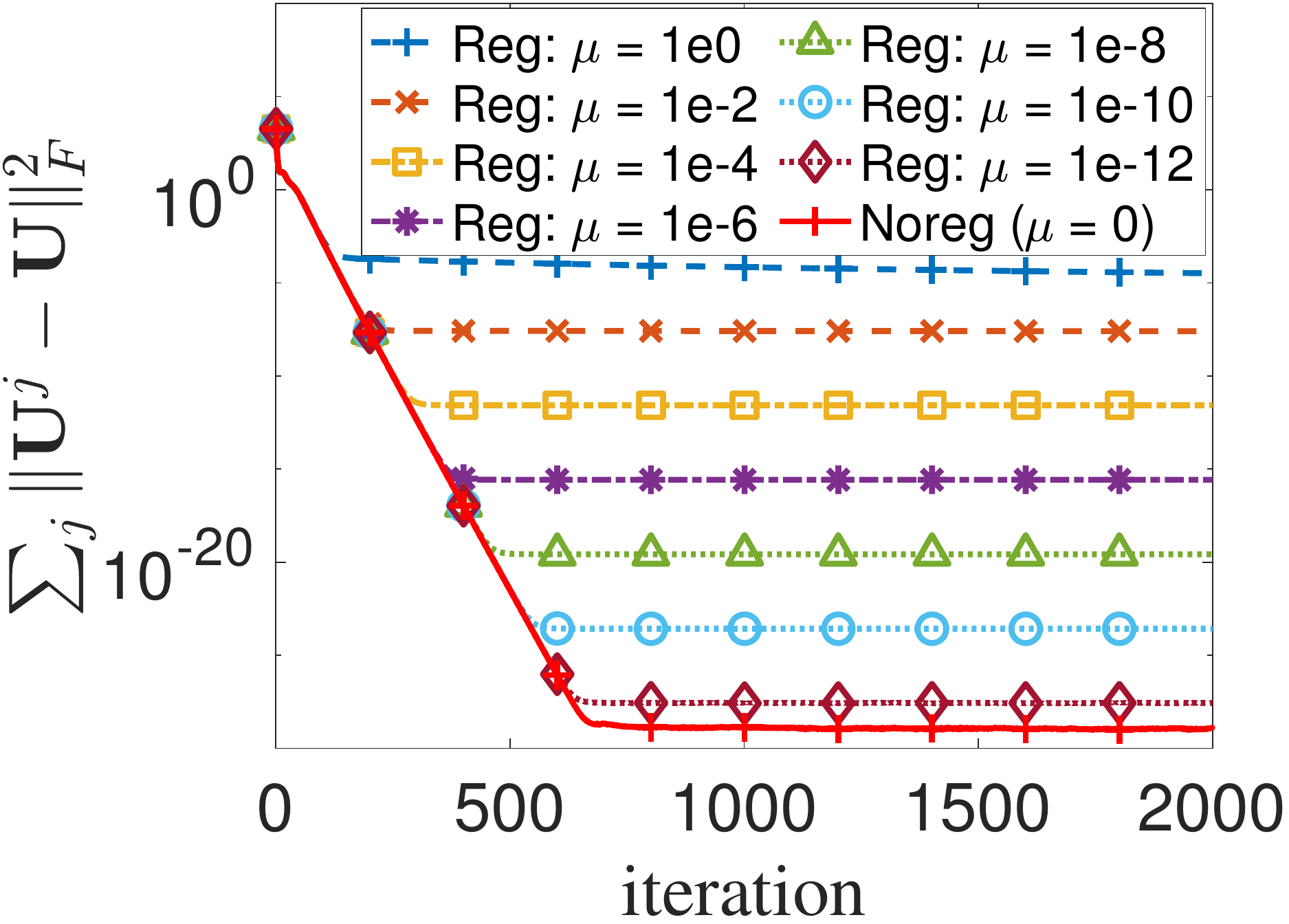}
\centerline{\footnotesize{(e) Distributed}}
\end{minipage}
\caption{Convergence of gradient descent and distributed gradient descent for solving the matrix sensing problems in terms of different optimality errors. Here, $m=50,~ n=40,~ r=5$, and the number of measurements $p=3\max\{m,n\}r$.}
\label{fig:ms}
\end{figure*}

\vspace{0.15cm}
\noindent{\bf Degenerate case: $\rank(\U\V^\top) < r$}

\vspace{0.15cm}


Note that
$
\rank(\mU^\top\mU\mV^\top\mV)\le\rank(\mU\mV^\top)<r,
$
which implies that
$
\det(\mU^\top\mU\mV^\top\mV)=\det(\mU^\top\mU)\det(\mV^\top\mV)=0.
$
Then, either $\det(\mU^\top\mU)=0$ or $\det(\mV^\top\mV)=0$. Or equivalently, either $\rank(\mU^\top\mU)<r$ or $\rank(\mV^\top\mV)<r$. Note that $\rank(\mU)=\rank(\mU^\top\mU)$ and $\rank(\mV)=\rank(\mV^\top\mV)$. Then, of the following two statements, \SLb{at least one of them is true.}
\begin{itemize}
\item[(i)] $\exists~\b\neq\zero~\text{such that}~\b\in\Null(\U)$, i.e., $\U\b = \zero$.
\item[(ii)] $\exists~\b\neq\zero~\text{such that}~\b\in\Null(\V)$, i.e., $\V\b = \zero$.
\end{itemize}

Note that for any critical point $(\U,\V)$, \MBW{either}
\begin{align*}
&	\nabla f(\mU\mV^\top)=\zero \Rightarrow (\U,\V) \tmop{is~a~global~minimum}, \tmop{\MBW{~or}}
\\
& \nabla f(\mU\mV^\top)\neq\zero \Rightarrow \exists ~(i,j),~~  \lg \nabla f(\mU\mV^\top),\e_i\e_j^\top \rg \neq 0.
\end{align*}
Next, we focus on the second case and show that such kinds of critical points are strict saddles.

Assume that (i) is true. Construct $\D = [\D_{\U}^\top~\D_{\V}^\top]^\top$ with $\D_{\U}^\top = \b \e_i^\top \in \R^{r\times n}$ and
$\D_{\V}^\top = (\alpha\b) \e_j^\top \in \R^{r \times m}$.
Then, we have
$
\D_{\U}\D_{\V}^\top	= \alpha \|\b\|_2^2 \e_i \e_j^\top,
$ and
$\U \D_{\V}^\top =\U (\alpha \b)\e_j^\top = \zero$.
Plugging into the bilinear form of the Hessian, we get
{\small
\begin{align*}
&\nabla^2 g(\U,\V)\left[\D,\D\right]
\nn\\
=& \lg 2\nabla f(\mU\mV^\top), \mD_\mU\mD_\mV^\top\rg+  \nabla^2f(\mU\mV^\top)\left[\mPi,\mPi\right]
\nn\\
= & \lg 2\nabla f(\mU\mV^\top), \alpha \|\b\|_2^2 \e_i \e_j^\top \rg \!+\!  \nabla^2f(\mU\mV^\top)\!\left[\e_i(\mV\b)^\top,\e_i(\mV\b)^\top\right]
\nn\\
= & 2\alpha \|\b\|_2^2 \lg \nabla f(\mU\mV^\top),  \e_i \e_j^\top \rg \!+\!  \nabla^2f(\mU\mV^\top)\!\left[\e_i(\mV\b)^\top,\e_i(\mV\b)^\top\right]
\end{align*}
}
\hspace{-.3cm}
Now using the fact that $ \lg \nabla f(\mU\mV^\top\SLb{)}, \e_{i}\e_j^\top\rg \neq 0$, $\|\b\|_2^2\neq 0$ and that $\nabla^2f(\mU\mV^\top)\left[\e_i(\mV\b)^\top,\e_i(\mV\b)^\top\right]$ is constant with respect to $\alpha$, we can always choose $\alpha$ in order to let the first term $\alpha \|\b\|_2^2 \lg 2\nabla f(\mU\mV^\top),  \e_i \e_j^\top \rg$  be negative enough so that $\nabla^2 g(\U,\V)\left[\D,\D\right]$ is negative. Therefore, we can conclude that such a critical point $(\mU,\mV)$ is a strict saddle.
Similarly, we can consider the case when (ii) is true and finish the proof.
\end{proof}


\vspace{-.05in}
\MBW{We note that, in Theorem~\ref{thm:local_mtxsen}, the requirement $\beta/\alpha\le 3/2$ is the same as in \cite{zhu2018global} where matrix sensing and other low-rank matrix recovery problems are discussed.}

\subsection{Landscape of Distributed Low-rank Matrix Recovery}
\label{sec:land_dis}

The following corollary extends the benign geometry to the {\em distributed} setting introduced in Section~\ref{sec:prob}. 
\begin{Corollary} \label{corr:distributed}	
\MBW{Under the assumptions in \Cref{thm:local_mtxsen}}, any critical point $(\{\mU^j\},\{\mV_j\})$ of the distributed problem  \eqref{eq:cost_nonconv_dis0:ECD:Penalty} satisfies  $\mU^{1}=\mU^{2}=\cdots=\mU^{J}=\mU$ for some $\mU$, and $(\mU,\{\mV_j\})$ is either a strict saddle or a global minimizer of~\eqref{eq:cost_nonconv_dis0:ECD:Penalty}.
\end{Corollary}

This result follows from Theorem~\ref{thm:local_mtxsen} and Lemma~\ref{thm:penalty:consensus}, which is a summary of~\cite[Proposition~2.3 and Theorem~2.7]{zhu2019distributed}.

\begin{Lemma}\cite{zhu2019distributed}
Let $\{w_{i,j}\}_{(i,j)\in\calG}$ be symmetric positive weights on any connected network $\calG$ over $[J]^2$. Then, {\em (i)} any critical point $(\{\mU^j\},\{\mV_j\})$ of the distributed problem \eqref{eq:cost_nonconv_dis0:ECD:Penalty} satisfies  $\mU^{1}=\mU^{2}=\cdots=\mU^{J}=\mU$ for some $\mU$, and {\em (ii)} $(\mU,\{\mV_j\})$ is a critical point of the centralized problem \eqref{eq:cost_nonconv_dis00}.
\label{thm:penalty:consensus}
\end{Lemma}

%
%
%
%
%



\section{Simulation Results}
\label{sec:simu}


In this section, we conduct several experiments 
to further support our theory. In particular, we first consider the following {\em centralized} matrix sensing problem
\begin{align}
\minimize_{\U\in\R^{n\times r},\V\in\R^{m\times r}}~g(\mU,\mV)\doteq \frac{1}{2}\|\calA(\mU\mV^\top-\mX^\star)\|_2^2,	
\label{eq:ms_cost}
\end{align}
where $\calA:\R^{n\times m}\to \R^p$ is a linear sensing operator, and $\mX^\star\in\R^{n\times m}$ is the true low-rank matrix with $\rank(\mX^\star)=r$. In order to compare the non-regularized setting with the regularized setting, we apply gradient descent with random initialization to minimize the following regularized cost function
\begin{align*}
\widetilde{g}(\mU,\mV)	=  \frac{1}{2}\|\calA(\mU\mV^\top-\mX^\star)\|_2^2	+\frac{\mu}{4}\|\mU^\top\mU-\mV^\top\mV\|_F^2
\end{align*}
with $\mu$ equal to $10^{-1},~10^{-2},~10^{-3},~10^{-4}$ and $0$. Note that the regularized cost function $\widetilde{g}(\mU,\mV)$ reduces to the  non-regularized cost function $g(\U,\V)$ when $\mu~=~0$. To set up the experiment, we choose $m=50,~ n=40,~ r=5$, and $p=3\max\{m,n\}r$. The true data matrix $\Xs$ is generated as $\Xs = \Us {\V^\star}^\top$, where $\Us$ and $\V^\star$ are two Gaussian random matrices \SLb{with entries following $\NN(0,1)$}. The linear sensing operator is generated as a $p \times nm$ Gaussian random matrix \SLb{with entries following $\NN(0,1)$}. We plot the fitting error $g(\mU,\mV)$ and the optimality error $\|\mU\mV^\top-\mX^\star\|_F^2$ as a function of the iteration number in \Cref{fig:ms} (a) and (b), respectively.
\MBW{One can observe that global optimality is achieved in all regularized and unregularized ($\mu=0$) cases.}\footnote{\MBW{All errors in the centralized experiments eventually decay below $10^{-20}$.}}
Therefore, in the case of {\em centralized} matrix sensing, the balancing regularizer $\mathrm{R}(\mU,\mV)$ in~\eqref{eq:balance_reg} is not necessary to obtain a benign landscape.

Next, we repeat the above experiments on the corresponding {\em distributed} matrix sensing problem, namely, minimizing the following cost functions
\begin{align*}
\widetilde{g}_d(\mU,\mV)	&\!=\!  g_d(\U,\V)	+\frac{\mu}{4}\|\mU^\top\mU-\mV^\top\mV\|_F^2,~\tmop{and}
\\
g_d(\mU,\mV)	&\!=\!  \frac 1 2\sumj \|\calA_j(\mU^j\mV_j^\top\!-\!\mX^\star_j)\|_2^2 +\!\!\!\!\!\sum_{(i,j)\in\calG}\!\!\!\!\! w_{i,j} \|\mU^j\!-\!\mU^i\|_F^2 .
\end{align*}
We set $\mu = 0$ and $10^{\alpha}$ with $\alpha = 0:-2:-12$. We choose $J=5$, $n=50$, $m_j = 5$, $m = \sumj m_j = 25$, $r = 5$, and $p=3\max\{m,n\}r$. We generate $\{w_{i,j}\}_{\calG}$ by performing hard thresholding on a random non-negative symmetric matrix \SLb{with off-diagonal entries being uniformly distributed random numbers in the interval $(0,1)$ and zero diagonal entries}.
Other parameters are set same as in the centralized framework.
We again present the fitting error $\sumj \|\calA_j(\mU^j\mV_j^\top-\mX^\star_j)\|_2^2$, the optimality error $\sumj\|\mU^j\mV_j^\top-\mX_j^\star\|_F^2$, and the {\em consensus} error $\sum_j \|\U^j-\U\|_F^2$ as a function of the iteration number in \Cref{fig:ms} (c), (d) and (e), respectively.
\MBW{While the existing literature does not guarantee a benign geometry for the distributed problem with regularization, we do see near-optimal convergence with near-consensus when the regularizer is sufficiently small. Our Corollary~\ref{corr:distributed} does apply to the distributed unregularized problem, and in this case we do indeed see global optimality and exact consensus.}

\section{Conclusion}
\label{sec:conc}

This work closes the theory-practice gap for the factorization approach in low-rank matrix optimization \SLb{when the cost function is restricted strongly convex and smooth} by showing that the balancing regularizer is not necessary in geometric analysis, in agreement with practical observations. We have proved that any critical point of the \MBW{unaltered factorized objective function (without regularizer)} is either a global minimum or a strict saddle in both {\em centralized} and {\em distributed} settings.


\newpage
\balance
\bibliographystyle{ieeetr}
\bibliography{SS_noreg}

\end{document}